\documentclass[aps,pre,reprint,a4paper,amsfonts,amsmath,amssymb,nofootinbib,10pt]{revtex4-2}
\usepackage[english]{babel}
\usepackage{graphicx}
\usepackage{mathrsfs}
\usepackage{enumerate}
\usepackage{hyperref}
\usepackage{etoolbox}

\begin{document}

\def\Journal#1#2#3#4#5#6#7{#1, #4 \textbf{#5}, #6 (#2).}
\def\JournalE#1#2#3#4#5#6{#1, #4 \textbf{#5}, #6 (#2).}
\def\Book#1#2#3#4#5{{#1}, {\it #3} (#4, #5, #2).}

\newcommand{\dd}{\mbox{d}}
\newcommand{\NN}{\mathbb{N}}
\newcommand{\ZZ}{\mathbb{Z}}
\newcommand{\RR}{\mathbb{R}}
\newcommand{\PP}{\mathbb{P}}
\newcommand{\EE}{\mathbb{E}}
\newcommand{\uu}{\mathbf{1}}
\newcommand{\RM}{\mathscr{R}}
\newcommand{\FM}{\mathscr{F}}
\newcommand{\NM}{\mathscr{N}}
\newcommand{\TT}{\mathscr{T}}
\newcommand{\step}{s}
\newcommand{\vv}{v_\text{d}}

\ifundef{\eqref}{\def\eqref#1{(\ref{#1})}}
\ifundef{\text}{\def\text#1{\textrm{#1}}}

\title{Merge of two oppositely biased Wiener processes}
\author{Miquel Montero}\email{miquel.montero@ub.edu}
\affiliation{Departament de F\'{\i}sica de la Mat\`eria Condensada, Universitat de Barcelona (UB), and Universitat de Barcelona Institute of Complex Systems (UBICS), Mart\'{\i} i Franqu\`es 1, E-08028 Barcelona, Spain}

\date{\today}

\begin{abstract}
We introduce a technique to merge two biased Brownian motions into a single regular process. The outcome follows a stochastic differential equation with a constant diffusion coefficient and a non-linear drift.  The emerging stochastic process has outstanding properties, such as spatial and temporal translational invariance of its mean squared displacement, and can be efficiently simulated via a random walk with site-dependent one-step transition probabilities.
\end{abstract}
\maketitle

\section{Introduction}

The time evolution of many physical systems can be expressed in terms of a stochastic differential equation (SDE) driven by the standard Wiener process \cite{G85,VK92,BO03}. In the case of the location of a particle, the equation may describe the behavior of its velocity under the action of both deterministic and random forces, when the inertial effects can be neglected. Then, the combined effect of the deterministic forces configures the \emph{drift} whereas the intensity of the random perturbations dictates the shape of the \emph{diffusion coefficient}. In the most general situation both magnitudes depend on the position of the particle~\cite{SSMKG82}, as in heterogenous media~\cite{CCM13,GD15,RDG17}; on the time variable, as when the system suffers of ageing~\cite{MJCB14}; or both, as in ecological problems~\cite{FM15,MGHT21}. 

However, in many practical situations, it is enough to assume that the diffusion coefficient is a constant that does not change either explicitly or implicitly in time, and that the drift is a function of the position of the particle, and even then obtain rich behaviour~\cite{CM95}. This is the case in which the Einstein equation applies: Particles moving by the action of an external force through a homogeneous viscous medium and experiencing drag, in such a way a terminal velocity should be reached. The fluctuation-dissipation theorem states that the strength of the fluctuations around this value is the product of the \emph{mobility} of the particles, that is, the ratio of the velocity to the applied force, times the absolute temperature of the medium, measured in energetic units~\cite{RK66}.

Even in this simplified scenario, with a constant diffusion coefficient and a position-dependent drift, the probability density function (PDF) of the location of the particle can be seldom found in a full explicit form, like in the Ornstein--Uhlenbeck process~\cite{OU30}. The equation that that governs the time evolution of the PDF associated to a SDE is the Fokker--Planck equation (FPE), a partial differential equation that frequently must be solved using numerical methods~\cite{R96,AS11,CF13}. Some notable exceptions, in addition to the more common and less informative situation of when the FPE admits a steady-state solution, include cases where the formal analogy between the FPE and the Schr\"{o}dinger equation can be exploited~\cite{BKM13,CCF14}, or a mixture model can be constructed~\cite{BM02,BMS03}.

Here we adopt an approach in line with the latter idea: the merge of two biased Brownian motions into a single process. The generated SDE has a constant diffusion coefficient with a hyperbolic tangent function as a drift. While one can find in the literature apparently similar processes~\cite{LC70}, they correspond to non-linear models for the restoring force associated with the stiffness of certain materials |see, e.g.,~\cite{YW21} and references therein| where the drift has a negative sign which results in a well different behavior: the process is then mean-reverting and admits a stationary probability distribution. In our case, in turn, the outcome shares properties with the random walk introduced in a previous work~\cite{MMa,MMb}.

The paper is structured as follows: In Sect.~\ref{Sec:Stochastic_process} we state the desired properties of the process (PDF, mean and variance) and find the SDE for the only process that satisfies these demands. In Sect.~\ref{Sec:MSD} we extend the statistical analysis of the process to the covariance and show how the mean squared displacement exhibits spatial and temporal translational invariance. Section~\ref{Sec:RW} introduces the connection between the continuous process and its discrete counterpart. The paper ends with Sect.~\ref{Sec:conclusions} where conclusions are drawn. 

\section{The stochastic process\label{Sec:Stochastic_process}}

Consider $X^{\pm}_t$, two biased Brownian motions starting at point $X^{\pm}_0=x_0$,
\begin{equation}
\dd X^{\pm}_t\equiv\pm \vv \dd t+ \sigma \dd W^{\pm}_t,
\label{eq:SDEpm}
\end{equation}
with $\vv$ and $\sigma$ positive constants, and $W^{\pm}_t$ two independent Wiener processes. The FPEs for the corresponding PDFs $p^{\pm}(x,t;x_0)$,
\begin{equation}
p^{\pm}(x,t;x_0)\dd x\equiv \PP\left(x<X^{\pm}_{t}\leq x+\dd x | X^{\pm}_0=x_0\right),
\end{equation}
the (conditional transition) probabilities $\PP\left(\cdot|\cdot\right)$, are
\begin{equation}
\frac{\partial}{\partial t} p^{\pm}(x,t;x_0)=\mp \vv \frac{\partial}{\partial x} p^{\pm}(x,t;x_0) +\frac{1}{2} \sigma^2 \frac{\partial^2}{ \partial x^2} p^{\pm}(x,t;x_0),
\label{eq:FPE_pm}
\end{equation} 
whose solutions read
\begin{equation}
p^{\pm}(x,t;x_0)=\frac{1}{\sqrt{2 \pi \sigma^2 t} }e^{-\frac{(x-x_0\mp \vv t)^2}{2 \sigma^2 t}},
\end{equation} 
that is, Gaussian distributions with the means and the variances equals to:
\begin{eqnarray}
\EE\left[X^{\pm}_t|X^{\pm}_{0}=x_0\right]=x_0\pm \vv t,\\
\EE\left[\left.\left(X_t^{\pm}\right)^2\right.|X^{\pm}_{0}=x_0\right]-\EE\left[X^{\pm}_t|X^{\pm}_{0}=x_0\right]^2=\sigma^2 t,
\end{eqnarray}
where we denote by $\EE\left[\cdot|\cdot\right]$ the conditional expectation of a random variable under the suitable probability measure.

Our objective is to construct a new stochastic process $X_t$, also starting at $X_0=x_0$, whose PDF $p(x,t;x_0)$,
\begin{equation}
p(x,t;x_0)\dd x\equiv \PP\left(x<X_{t}\leq x+\dd x | X_0=x_0\right),
\end{equation}
is
\begin{eqnarray}
p(x,t;x_0) &=& \frac{1+f(x_0)}{2}p^{+}(x,t;x_0) \nonumber\\
&+&\frac{1-f(x_0)}{2}p^{-}(x,t;x_0),
\label{eq:PDF_comp}
\end{eqnarray}
with $f(x_0)$ a smooth function of its argument that fulfills $-1\leq f(x_0) \leq 1$. A possible approach to the problem could be to define $X^{\Xi}_t$,
\begin{equation}
\dd X^{\Xi}_t \equiv  \Xi(x_0) \dd X^{+}_t +\left[1-\Xi(x_0)\right]\dd X^{-}_t ,
\label{eq:SDE_Bernoulli}
\end{equation}
where $\Xi(x_0)$ is a Bernoulli random variable independent of $W^{\pm}_t$, $\Xi(x_0)\in\{0,1\}$, with
\begin{equation}
\PP\left(\Xi(x_0)=1\right)= \frac{1+f(x_0)}{2}.
\end{equation}
The mean and variance of $X^{\Xi}_t$ are, respectively,
\begin{equation}
\EE\left[X^{\Xi}_t|x_0\right] = x_0+\vv f(x_0)  t,
\label{eq:mean_Benoulli}
\end{equation}
and 
\begin{equation}
\EE\left[\left(X^{\Xi}_t\right)^2|x_0\right]-\EE\left[X^{\Xi}_t|x_0\right]^2=\sigma^2 t+\left[1-f(x_0)^2\right] \vv^2 t^2,
\label{eq:variance_Bernoulli}
\end{equation} 
where we have simplified everywhere the notation relative to the condition, e.g., $\EE\left[X^{\Xi}_t|x_0\right]\equiv\EE\left[X_t^{\Xi}|X^{\Xi}_0=x_0\right]$. Hence, the time evolution of the variance of the process we are looking for will present two well-different regimes: The first term in Eq.~\eqref{eq:variance_Bernoulli} will be the most relevant initially but, as time increases, the second term will eventually dominate.

Equation \eqref{eq:SDE_Bernoulli} represents the random alternation of two independent processes as, e.g., the movement of (non interacting, identical) positively and negatively charged particles, $\pm q$, in a gas or liquid within a uniform electric field of magnitude $E$. In this case, in absence of noise, the moving particles will reach a terminal drift velocity $\pm \vv$ according to the formula
\begin{equation}
\vv =\mu_q E, 
\label{eq:drift_E}
\end{equation}
where $\mu_q$ is the electrical mobility, and $f(x_0)$ determines the relative density excess of one species with respect to the other. Therefore, $p(x,t;x_0)$ measures the joint probability of finding a particle, irrespective of its charge, at a given location of the medium. Moreover, if the origin of the noise is the fluctuation-dissipation relation on has the Einstein--Smoluchowski equation,
\begin{equation}
\sigma^2 =\frac{2 \mu_q k_\text{B} T}{q}, 
\label{eq:ES}
\end{equation}
where $k_\text{B}$ is the Boltzmann constant and $T$ is the absolute temperature.

Therefore, Eq.~\eqref{eq:SDE_Bernoulli} does not describe the evolution of a single system but a mixture of two of them, governed by the Bernoulli variable. A naive way to eliminate the dependence of the process in $\Xi(x_0)$ would be to replace the random variable by its mean value,
\begin{equation}
\dd X^*_t \equiv  \frac{1+f(x_0)}{2}\dd X^{+}_t +\frac{1-f(x_0)}{2}\dd X^{-}_t,
\label{eq:SDE_sup}
\end{equation}
but in this case $p^*(x,t;x_0)$ will not fulfill Eq.~\eqref{eq:PDF_comp}.
Indeed, $X^*_t$ is yet another biased Brownian motion, 
\begin{equation}
\dd X^*_t = \vv f(x_0) \dd t+ \sigma \sqrt{\frac{1+f^2(x_0)}{2}} \dd W^*_t,
\label{eq:SDE_sup_bis}
\end{equation}
with the right mean value and the wrong variance. In spite of that, Eq.~\eqref{eq:SDE_sup} is a source of inspiration for defining process $X_t$, 
\begin{eqnarray}
\dd X_t &\equiv& \frac{1+f(X_t)}{2}\left[\vv \dd t+ \sigma \dd W_t\right]\nonumber\\
&+&\frac{1-f(X_t)}{2}\left[-\vv \dd t+ \sigma \dd W_t\right]\nonumber\\
&=&\vv f(X_t) \dd t+ \sigma \dd W_t,
\label{eq:X_def}
\end{eqnarray}
where we remark that there is a single origin of Gaussian white noise, $W_t$. To ensure that Eq. \eqref{eq:PDF_comp} characterizes this process one needs that $f(X_t)$ is a \emph{martingale}, i.e.,
\begin{equation}
\EE\left[f(X_t)|x_0\right]=f(x_0),
\end{equation} 
since from Eqs.~\eqref{eq:mean_Benoulli} and~\eqref{eq:X_def} one has that
\begin{equation}
\dd \EE\left[X_t|x_0\right]=\vv \EE\left[f(X_t)|x_0\right] \dd t= \vv f(x_0) \dd t.
\end{equation}
Moreover, $f(x_0)$ must be bounded, $-1 \leq f(x_0)\leq 1$, otherwise the resulting $p(x,t;x_0)$ could become eventually negative for certain values of $x$. 

The stochastic differential equation that governs the evolution of $f(X_t)$ is:~\footnote{We have applied the It\^{o} formula, the same conclusion steams from the use of the Stratonovich calculus. Alternatively, one can resort to the Kolmogorov backward equation to obtain Eq.~\eqref{eq:ODE}.}   
\begin{eqnarray}
\dd f(X_t) &=& f'(X_t) \dd X_t + \frac{1}{2} \sigma^2 f''(X_t) \dd t\nonumber\\
&=&\left[\vv f'(X_t) f(X_t) +\frac{1}{2} \sigma^2 f''(X_t)\right] \dd t\nonumber\\
&+& \sigma f'(X_t) \dd W_t.
\end{eqnarray}
Then, $f(X_t)$ is a martingale if and only if the term inside the square brackets is null, that is, if $f(u)$ is the solution of the following ordinary differential equation 
\begin{equation}
f''(u)+2 f'(u) f(u)=0,
\label{eq:ODE}
\end{equation} 
where we have defined the dimensionless variable $u\equiv \kappa x$ in terms of $\kappa$,
\begin{equation}
\kappa \equiv \frac{\vv}{\sigma^2}.
\label{eq:kappa_def}
\end{equation} 
The general solution of Eq. \eqref{eq:ODE} is
\begin{equation}
f(u)=k_0\frac{\sinh(k_0 u)+k_1\cosh(k_0 u)}{\cosh(k_0 u)+k_1\sinh(k_0 u)}
\label{eq:ODE_sol_gen}
\end{equation} 
with $k_0$ and $k_1$ arbitrary complex constants, $k_0=a+i\alpha$, $k_1=b+i\beta$ . In order to restrict $f(u)$ to the reals, one must have that $k_0$ and $k_1$ are either both real or both imaginary, i.e., either $\alpha=\beta=0$, or $a=b=0$. In the former case one has
\begin{equation}
f(u)=\frac{\sinh(u)+b\cosh(u)}{\cosh(u)+b\sinh(u)},
\label{eq:ODE_sol_real}
\end{equation} 
where we have set $a=1$, since this numerical constant can be absorbed in the definition of $\vv$. After that, the functional form of Eq.~\eqref{eq:ODE_sol_real} depends on the value of $b$. For $b= \pm 1$ one recovers the original biased processes with velocity $\pm \vv$.  For $|b|<1$ one gets
\begin{equation}
f(u)=\tanh(u+u_b),
\label{eq:ODE_sol_tanh}
\end{equation} 
with $b=\tanh(u_b)$, while for $|b|>1$ one has
\begin{equation}
f(u)=\coth(u+u_b),
\label{eq:ODE_sol_coth}
\end{equation} 
with $b=\coth(u_b)$. This latter solution is precluded within our framework because we must have $|f(u)|\leq 1$ if we want that Eq.~\eqref{eq:PDF_comp} is well defined and, therefore, its is left for future investigations. Something similar happens in the complementary case, $a=b=0$, where Eq.~\eqref{eq:ODE_sol_gen} reduces to
\begin{equation}
f(u)=- \tan(u+u_\beta),
\label{eq:ODE_sol_tan}
\end{equation} 
once we have set $\alpha=1$ and $\beta=\tan(u_\beta)$. Here we can force the constraint $|f(u)|\leq 1$ by restricting the domain of $u$, for example with $|u|\leq \pi/4$ and $\beta=0$. However, this renders invalid Eq.~\eqref{eq:PDF_comp} as well, since  $p^{\pm}(x,t;x_0)$ have support in the whole real line. Indeed, this choice for $f(u)$ is related with a previous work, see~\cite{MMb}.

In conclusion, only for~\footnote{Since $p^{\pm}(x,t;x_0)$ are translational invariant, we can freely choose $b=0$. This is one of the reasons why it has not been set $x_0=0$. } 
\begin{eqnarray}
\dd X_t &=&\vv \tanh\left(\kappa X_t\right) \dd t+ \sigma \dd W_t,
\label{eq:SDE_final}
\end{eqnarray}
one has that $p(x,t;x_0)$ is positively definite, and shows the desired structure,
\begin{eqnarray}
p(x,t;x_0)&=&  \frac{e^{\kappa x_0}}{e^{\kappa x_0}+e^{-\kappa x_0}}p^{+}(x,t;x_0) \nonumber\\
&+& \frac{e^{-\kappa x_0}}{e^{\kappa x_0}+e^{-\kappa x_0}} p^{-}(x,t;x_0) \nonumber\\
&=&\frac{1}{\sqrt{2 \pi \sigma^2 t} } \frac{ \cosh (\kappa x)}{\cosh (\kappa x_0)}  e^{-\frac{(x-x_0)^2}{2 \sigma^2 t}-\frac{\vv^2 t}{2 \sigma^2} },
\label{eq:PDF_sol}
\end{eqnarray}
with the appropriate mean
\begin{equation}
\EE\left[X_t|x_0\right]=  x_0+\vv \tanh(\kappa x_0)  t,
\label{eq:mean}
\end{equation}
and variance
\begin{eqnarray}
\text{var}\left[X_t|x_0\right]&\equiv& \EE\left[X_t^2|x_0\right]-\EE\left[X_t|x_0\right]^2\nonumber\\
&=&\sigma^2 t+\left[\frac{\vv t}{\cosh(\kappa x_0)}\right]^2.
\label{eq:variance}
\end{eqnarray} 
(Recall that $\vv =\kappa \sigma^2$: we will use together $\vv$, $\kappa$ and $\sigma$ in the same equation whenever this clarifies the expression.) It can be checked by direct insertion how Eq.~\eqref{eq:PDF_sol} satisfies the FPE associated to Eq.~\eqref{eq:SDE_final}:
\begin{eqnarray}
\frac{\partial}{\partial t} p(x,t;x_0)&=&-\vv  \frac{\partial}{\partial x}\left[\tanh(\kappa x) p(x,t;x_0)\right] \nonumber\\
&+&\frac{1}{2}\sigma^2  \frac{\partial^2}{ \partial x^2} p(x,t;x_0).
\label{eq:FPE_final}
\end{eqnarray}
Incidentally, since the drift and the diffusion coefficient are Lipschitz continuous functions, the solution of Eq.~\eqref{eq:SDE_final}, given $x_0$, is continuous and unique~\cite{BO03}.

If we continue with the example in which $X_t$ describes the position of charged particles moving within a gas or liquid, now all the particles are identical, with charge $q$, and $x=0$ marks the point where there is a sudden change in the direction of the electric field, that passes from $-E$, for $x<0$, to $E$, for $x>0$, with the sign prescription that makes $q E>0$: note that $\kappa$ is positive by definition, cf. Eq.~\eqref{eq:kappa_def}, and here
\begin{equation}
\kappa = \frac{q E}{2 k_\text{B}T}.
\label{eq:kappa_ex}
\end{equation}

\section{Mean squared displacement\label{Sec:MSD}}

Let us progress with the analysis of the expected values of $X_t$ and consider now the covariance,
\begin{equation}
\text{cov}\left[X_{t+\tau}X_{t}|x_0\right]\equiv \EE\left[X_{t+\tau}X_{t}|x_0\right]-\EE\left[X_{t+\tau}|x_0\right]\EE\left[X_{t}|x_0\right],
\label{eq:cov_def}
\end{equation}
with $\tau\geq 0$. The process is time homogeneous and one can easily check that the PDF in \eqref{eq:PDF_sol} satisfies the corresponding Chapman--Kolmogorov equation~\cite{VK92}:
\begin{equation}
p(x,t;x_0)=\int_{-\infty}^{+\infty} p(x,t-t';x') p(x',t';x_0) \dd x', 
\end{equation} 
for any $0\leq t' \leq t$. Therefore, to compute Eq.~\eqref{eq:cov_def} one can use the tower property of the conditional expectation and get
\begin{eqnarray}
\EE\left[X_{t+\tau}X_{t}|x_0\right]&=&\EE\left[\left. \EE\left[X_{t+\tau}X_{t}|X_{t}\right]\right|x_0\right]\nonumber \\
&=&\EE\left[\left. X_t \left(X_{t}+\vv \tanh(X_{t})\tau\right) \right|x_0\right]\nonumber\\
&=&x_0^2+\left[\sigma^2 + 2x_0 \vv \tanh(\kappa x_0)\right]t +\vv^2 t^2 \nonumber \\ 
&+&\vv \tau  \left[x_0 \tanh(\kappa x_0) + \vv t \right],
\label{eq:corr}
\end{eqnarray}
and, from this,
\begin{equation}
\text{cov}\left[X_{t+\tau}X_{t}|x_0\right]=\sigma^2 t + \left[\frac{\vv}{\cosh(\kappa x_0)}\right]^2 t(t+\tau),
\label{eq:cov}
\end{equation}
which closely resembles Eq.~\eqref{eq:variance}. If one defines the squared displacement of the process between times $t$ and $t+\tau$, $\Delta X^2(t,t+\tau)$, as
\begin{equation}
\Delta X^2(t,t+\tau)\equiv  (X_{t+\tau}-X_{t})^2,
\label{eq:MSD}
\end{equation}
one can use Eq. \eqref{eq:corr} to obtain the mean squared displacement (MSD) of the process,
\begin{eqnarray}
\EE\left[\Delta X^2(t,t+\tau)|x_0\right]
=\sigma^2 \tau + \vv^2 \tau^2,
\label{eq:EMSD}
\end{eqnarray}
independent of $t$ and $x_0$. It can easily be verified that this result is the one which would be obtained from any of the biased Brownian motions, a distinctive trait~\cite{DNNPR01}. 

Indeed, if one computes the MSD for a general process satisfying Eq. \eqref{eq:X_def} one gets
\begin{eqnarray}
\EE\left[\Delta X^2(t,t+\tau)|x_0\right]= \sigma^2 \tau \nonumber\\
+ 2\vv \left(\int_{t}^{t+\tau} \EE\left[X_{t'}f(X_{t'})|x_0\right] \dd t'-  \tau \EE\left[X_{t}f(X_{t})|x_0\right] \right), \nonumber\\
\label{eq:EMSD_gen}
\end{eqnarray}
which will be both independent of $t$ and $x_0$ if and only if
\begin{equation}
\EE\left[X_{t}f(X_{t})|x_0\right]=x_0 f(x_0) + c_0 \vv t,
\label{eq:EXfX}
\end{equation}
where $c_0$ is a numerical constant. 
From the Feynman--Kac formula, one has that $\EE\left[X_{t}f(X_{t})|x_0\right]$ must fulfill the Kolmogorov backward equation, that for our time-homogeneous process reads:
\begin{eqnarray}
\frac{\partial}{\partial t} \EE\left[X_{t}f(X_{t})|x_0\right]&=&\vv f(x_0)\frac{\partial}{\partial x_0} \EE\left[X_{t}f(X_{t})|x_0\right]\nonumber\\
&+&\frac{1}{2}\sigma^2  \frac{\partial^2}{ \partial x_0^2} \EE\left[X_{t}f(X_{t})|x_0\right],
\label{eq:KBE}
\end{eqnarray}
which leads to
\begin{equation}
c_0\vv=\vv f(x_0) \frac{\dd}{\dd x_0} \left[x_0 f(x_0)\right]+ \frac{1}{2} \sigma^2 \frac{\dd^2}{ \dd x_0^2} \left[x_0f(x_0)\right],
\end{equation}
that is,
\begin{equation}
u g'(u)+2 g(u) = 2 c_0,
\label{eq:g}
\end{equation}
with $u=\kappa x_0$ in this case, and 
\begin{equation}
g(u)\equiv f'(u)+f(u)^2.
\end{equation}
The solution of Eq. \eqref{eq:g} is equal to
\begin{equation}
g(u) = c_0+\frac{c_1}{u^2},
\label{eq:g_sol}
\end{equation}
where $c_1$ is a new integration constant. With $c_1=0$ we recover the cases we have already analyzed, cf. Eq.~\eqref{eq:ODE}, while $c_1 \neq 0$, even when $c_0=0$, leads to new scenarios that are beyond the scope of the present analysis but deserve future attention: the search for regularities in the properties of the MSD is a topic of current active interest~\cite{AG12, MJCB14}.

\section{Simulation\label{Sec:RW}}

In order to efficiently simulate possible realizations of $X_t$ one can resort to standard techniques~\cite{KP92,TS12}, as the Euler--Maruyama~\cite{GM55} method or the stochastic Runge--Kutta scheme~\cite{RH92}. However, the process can be dealt as the limit of the one-dimensional random walk (RW) we have analyzed in the past \cite{MMa,MMb}. There, we introduced an infinite Markov chain, whose one-step evolution can be expressed as follows: If at time $t$ the walker is at a given location, $X_{t}=x$, then at time $t+\Delta t$ one has
\begin{equation}
X_{t+\Delta t}=\left\{
\begin{array}{ll}
x+\Delta x,&\mbox{ with probability } p_{x\to x+\Delta x},\\
x-\Delta x,&\mbox{ with probability }  p_{x\to x- \Delta x},
\end{array}
\right.
\label{process}
\end{equation}
with the following inhomogeneous, one-step transition probabilities:
\begin{equation}
p_{x \to x\pm \Delta x}=\frac{1}{2} \frac{\cosh((x\pm\Delta x)\xi/\Delta x)}{\cosh (\xi ) \cosh (\xi x/\Delta x)},
\label{eq:pn}
\end{equation} 
where $\xi >0$ is a parameter that controls all the transition probabilities. (The case $\xi=0$ leads to the standard RW.) This choice for $p_{x \to x\pm \Delta x}$ defines valid transition probabilities, i.e., they are positive and the total probability leaving any given site equals to $1$, 
\begin{equation}
p_{x \to x+ \Delta x}+p_{x \to x- \Delta x}=1,
\label{eq:conservation}
\end{equation} 
with the remarkable property that the probability of performing a closed loop is independent of $x$: 
\begin{equation}
p_{x \to x\pm \Delta x} \cdot p_{x\pm \Delta x\to x}=\frac{1}{4 \cosh^2 (\xi )},
\label{eq:one_loop}
\end{equation} 
for one-step loops, which simply translates to loops of any size and topology. In the original references, for a matter of simplicity, it was assumed that $\Delta t=\Delta x =1$,  $X_t\in\ZZ$ and $t\in \NN_0$, but if one keeps these magnitudes free, it can be shown that one has
\begin{eqnarray}
\PP\left( X_{t}=x| X_0=x_0\right)\nonumber \\
\equiv{\frac{t}{\Delta t} \choose \frac{t}{2 \Delta t} -\frac{x-x_0}{2 \Delta x}} \frac{ \cosh (\xi x/\Delta x)}{\left[2 \cosh (\xi )\right]^{t/\Delta t}  \cosh (\xi x_0/\Delta x)},
\label{eq:tp_gen}
\end{eqnarray} 
for general values of $x_0$, $x$ and $t$, as long as $t/\Delta t$ and $t/\Delta t-|x-x_0|/\Delta x$ are non-negative even integers. In the limit where $\Delta x\to 0$, $\Delta t\to 0$, and $\xi \to 0$, with   
\begin{equation*}
\sigma\equiv \frac{\Delta x}{\sqrt{\Delta t}}, \text{and} \; \kappa\equiv \frac{\xi}{\Delta x}
\end{equation*}
constant magnitudes, one recovers Eq.~\eqref{eq:PDF_sol}: the normal distribution appears by using the familiar Stirling approximation on the \emph{pure combinatorial} part of Eq.~\eqref{eq:tp_gen}, that is, the one that does not depend on  $\xi$, 
\begin{equation*}
{\frac{t}{\Delta t} \choose \frac{t}{2 \Delta t} -\frac{x-x_0}{2 \Delta x}} \frac{1}{2^{t/\Delta t}}  \to \frac{1}{\sqrt{2 \pi \sigma^2 t} } e^{-\frac{(x-x_0)^2}{2 \sigma^2 t}} \dd x,
\end{equation*}
where we have made the identification $ 2\Delta x \to \dd x$ since this corresponds to the minimum change in $x$, while
\begin{equation*}
\frac{ \cosh (\xi x/\Delta x)}{\left[\cosh (\xi )\right]^{t/\Delta t}  \cosh (\xi x_0/\Delta x)} \to \frac{ \cosh (\kappa x)}{\cosh (\kappa x_0)} e^{-\frac{1}{2} \kappa^2 \sigma^2 t}.
\end{equation*}
This refers to the PDF of the location of the process but the analogy is not restricted to it, since the attributes of this RW reported in \cite{MMa} have their counterpart in the continuous process $X_t$, as the alluded spatial and temporal translational invariance of the MSD or the connection between Eq.~\eqref{eq:PDF_sol} and the first-passage time density. We leave for a future publication a detailed account of theses properties.

\section{Conclusions\label{Sec:conclusions}}

We have introduced a technique to merge two Wiener processes with opposite biases. Albeit the adopted approach can produce additional models with similar properties, we have proven that the stochastic process with a hyperbolic tangent as a drift is the only one that satisfies all the demands. Moreover, it belongs to the larger set of processes with a constant diffusion coefficient whose mean squared displacement is a function exclusively of the time lag. We have mathematically characterized this family, a family that has yet to be explored in depth.

Finally, the proposed technique can be useful to address similar problems where two related but different dynamics are wanted to be joined in a single process.

\acknowledgments
This research was funded by MCIN (Spain), Agencia Estatal de Investigaci\'on (AEI), grant number PID2019-106811GB-C33 (AEI/10.13039/501100011033); and by Generalitat de Catalunya, Ag\`encia de Gesti\'o d'Ajuts Universitaris i de Recerca (AGAUR), grant number No. 2017 SGR 1064.

\end{document}